\documentclass[12pt]{article}

\usepackage{amssymb}
\usepackage{epsfig}
\usepackage{color}

\newcommand{\ud}{\mathrm{d}}
\newcommand{\bean}{\begin{eqnarray*}}
\newcommand{\eean}{\end{eqnarray*}}
\newcommand{\bea}{\begin{eqnarray}}
\newcommand{\eea}{\end{eqnarray}}
\newcommand{\bbE}{{\mathbb{E}}}
\newcommand{\Exp}{E}
\newcommand{\bbP}{{\mathbb{P}}}
\newcommand{\Var}{{\mathrm{Var}}}
\def\Pr{P}
\newcommand{\R}{\mathbb{R}}
\newcommand{\Z}{\mathbb{Z}}
\newcommand{\F}{\mathcal{F}}
\newcommand{\N}{\mathbb{N}}
\newcommand{\eps}{{\varepsilon}}
\newcommand{\eqd}{\stackrel{\mathcal{D}}{=}}
\newtheorem{theorem}{Theorem}

\newtheorem{lemma}{Lemma}

\newtheorem{definition}{Definition}
\newcommand{\remark}{\noindent \textbf{Remark. }}
\newcommand{\remarks}{\noindent \textbf{Remarks. }}
\newcommand{\proof}{\noindent \textbf{Proof. }}
\newcommand{\sign}{\textrm{sign}}

\begin{document}

\title{Random walk in random environment with asymptotically zero perturbation}
\date{July 2006}

\author{M.V.~Menshikov\footnote{e-mail: \texttt{mikhail.menshikov@durham.ac.uk}}
 \\
 \normalsize
 Department of Mathematical Sciences,
 University of Durham,\\
\normalsize
 South Road, Durham DH1 3LE, England. 
\and Andrew R.~Wade\footnote{e-mail: \texttt{andrew.wade@bris.ac.uk}}
\\
\normalsize
 Department of Mathematics,
 University of Bristol,\\
\normalsize
 University Walk, Bristol BS8 1TW, England.}
 
\date{May 2006}

\maketitle

\begin{abstract}
We give criteria for ergodicity, transience and null recurrence
for the random walk in random environment on
$\Z^+=\{0,1,2,\ldots\}$, with reflection at the origin, where the
random environment is subject to a vanishing perturbation. Our
results complement existing criteria for random walks in random
environments and for Markov chains with asymptotically zero drift,
and are significantly different to these previously studied cases.
Our method is based on a martingale technique --- the method of
Lyapunov functions.
\end{abstract}

\vskip 3mm


\vskip 3mm


\section{Introduction}

In this paper we study a problem with a classical flavour that lies
in the intersection of two well-studied problems, those of random walks in one-dimensional
random environments
and Markov chains with asymptotically small drifts. Separately, these two problems
have received considerable attention,
but the problem considered in this paper has not
been analysed before. Further,
our results show that the system studied here
exhibits behaviour that is significantly different to
that of those previously studied systems.

The random walk in random environment (or RWRE for short) was
first studied by Kozlov \cite{kozlov} and Solomon \cite{solomon},
and has since received extensive attention; see for example
\cite{rev} or \cite{zeit} for surveys. This paper analyses the
behaviour of the RWRE for which the random environment is
perturbed by a vanishingly small amount.

The analysis of zero drift random walks in
two or more dimensions
by the method of Lyapunov functions
demonstrated the importance of the investigation of
one-dimensional stochastic
processes with asymptotically small drifts (see for example \cite{aim,lamp1,mai}). For example, if $(Z_t)$, with
$t=0,1,2,3\ldots$ time
is a random walk (with zero drift) in the nonnegative quarter plane, analysis
of the stochastic process $\|Z_t\|$, where $\|\cdot\|$ denotes
the Euclidean norm, involves the study of stochastic processes on the half-line
with mean drift asymptotically zero. 

Early work in this field was done by Lamperti \cite{lamp1,lamp2}.
Criteria for recurrence and transience are given in \cite{mai}, where the behaviour
in the critical
regime that Lamperti did not cover was also analysed.
Passage-time
moments are considered in \cite{aim}. In much of this work, Lyapunov functions
play a central role.

In this paper we demonstrate the essential difference between a
nearest-neighbour random walk in a deterministic environment,
perturbed from its critical (null-recurrent) regime, and a
nearest-neighbour random walk in a random environment, also
perturbed from its critical regime (sometimes called Sinai's
regime -- see below). Our results quantify the fact that in some
sense the random environment is more stable, in that a much larger
perturbation is required to disturb the null-recurrent situation.
In particular, we give criteria for ergodicity (i.e.~positive
recurrence), transience and null-recurrence for our perturbed
random walk in random environment. We will show that in our
(random environment) case the critical magnitude for the
perturbation is of the order of $n^{-1/2}$ (see Theorem
\ref{cor2}), where $n$ is the distance from the origin (in fact,
our more general results are much more precise than this). This
compares to a critical magnitude of the order of $n^{-1}$ in the
non-random environment case (see \cite{mai}, and Theorem
\ref{0928t} below).

Our method is based upon the theory
of Lyapunov functions, a powerful tool in the
classification of countable
Markov chains (see \cite{fmm}).
Such methods have proven effective
in the analysis of random walks
in random environments (see e.g.~\cite{cmp}),
in addition to Markov chains in non-random
environments.

Loosely speaking, motivation for our model
comes from some one-dimensional physical systems, such as
a particle
performing a random walk in a homogeneous random one-dimensional
field, subject to some vanishing perturbation (such as the
presence of another particle). Under what conditions is the perturbation
sufficient to alter the character of the random walk?

We now introduce the probabilistic model that we consider. First,
we need some notation. We introduce the function $\chi$ as follows, which
determines our perturbation as described below. Let $\chi : [0,\infty) \to [0,\infty)$
 be a function such that
\bea
\label{chicon}
\lim_{x \to \infty} \chi(x) =0.
\eea
As we shall see below, the property (\ref{chicon})
means that our perturbation is asymptotically small.

Here, we are interested in the one-dimensional
RWRE on the nonnegative integers (we use the notation
$\Z^+ := \{0,1,2,\ldots\}$), with reflection at the origin.
One can readily obtain results for the one-dimensional
RWRE on whole of $\Z$ in a similar manner. Formally, we define
our RWRE as follows.

We define sequences of random variables
$\xi_i$, $i=1,2,\ldots$ and $Y_i$, $i=1,2,\ldots$,
on some probability space $(\Omega,\F,\bbP)$, with the following properties.

Fix $\eps$ such that
$0<\eps<1/2$. Let $\xi_i$, $i=1,2,\ldots$, be a sequence of
i.i.d.~random
variables such that
\bea
\label{ue}
 \bbP [ \eps \leq \xi_1 \leq 1-\eps ] =1. \eea
The condition (\ref{ue}) is
sometimes
referred to as {\em uniform ellipticity}.

Let $Y_i$, $i=1,2,\ldots$, be another sequence of i.i.d.~random
variables,
taking values in $[-1,1]$,
on the same probability space as the $\xi_i$. We allow $Y_i$ to depend on $\xi_i$,
but any collections $(Y_{i_1},Y_{i_2}, \ldots, Y_{i_k})$,
$(\xi_{j_1},\xi_{j_2}, \ldots, \xi_{j_{k'}})$ are independent
if $\{ i_1, \ldots , i_k \} \cap \{ j_1, \ldots, j_{k'} \} =
\emptyset$.

For a particular realization of the
sequences $(\xi_i; i=1,2,\ldots)$ and $( Y_i; i=1,2,\ldots)$,
we define the quantities $p_n$ and $q_n$, $n=1,2,3,\ldots$ as follows:
\bea
\label{1006b}
 p_n & := & \left\{ \begin{array}{ll}
  \xi_n +Y_n\chi(n)  & ~~~~{\rm if}~~~ \eps/2 \leq \xi_n +Y_n \chi(n) \leq 1-(\eps/2) \\
  \eps/2  & ~~~~{\rm if}~~~ \xi_n +Y_n \chi(n) < \eps/2 \\
  1-(\eps/2)  & ~~~~{\rm if}~~~ \xi_n +Y_n \chi(n) > 1-(\eps/2) \end{array} \right.
 \nonumber\\
 q_n & := & 1-p_n.
\eea
We call a particular realization of $(p_n,q_n)$, $n=1,2,\ldots$, our
{\em environment}, and we denote it by $\omega$.
A given $\omega$ is then a realization of our random environment,
and is given in terms of the $\xi_i$ and $Y_i$ as in (\ref{1006b}). 

For a given environment $\omega$, that is, a realization of $(p_n,q_n)$, $n=1,2,\ldots$,
we
define
 the Markov chain $(\eta_t(\omega);t \in \Z^+)$ on
$\Z^+$, starting at some point in $\Z^+$, defined as follows:
$\eta_0(\omega)=r$ for some $r \in \Z^+$, and for $n=1,2,\ldots$,
\bea \label{1006bb}
 \Pr [ \eta_{t+1}(\omega) = n-1 | \eta_t(\omega) = n] & = & p_n,\nonumber\\
 \Pr [ \eta_{t+1}(\omega) = n+1 | \eta_t(\omega) = n] & = & q_n,
\eea
 and $\Pr [ \eta_{t+1}(\omega)=1 | \eta_t(\omega)=0] = 1/2$, $\Pr [ \eta_{t+1}(\omega)=0 |
\eta_t(\omega)=0] =1/2$. (Here $\Pr$ is the so-called {\em quenched} probability measure, 
i.e.~for a fixed environment $\omega$.)
The given
 form for the reflection at the origin
ensures that the Markov chain is {\em aperiodic}, which eases some technical complications, but
this choice is not special; it can be changed without affecting our results.

Recall that, from (\ref{chicon}), $\chi(n) \to 0$ as $n \to
\infty$. Thus, there exists $n_0 \in (0,\infty)$ such that, for
all $n \geq n_0$, $\chi(n) < \eps/2$. Hence, under condition
(\ref{ue}),
 for $\bbP$-almost every $\omega$ we have
$(\eps/2) < \xi_n +Y_n \chi(n) < 1-(\eps/2)$ (since the $Y_n$ are bounded). (For the remainder
of the paper we often use `a.e.~$\omega$' as shorthand for
`$\bbP$-almost every $\omega$' when the context is clear.)
Thus,
 for all $n \geq n_0$, (\ref{1006b}) implies that, for a.e.~$\omega$,
\bea
\label{1006c}
p_n = \xi_n +Y_n\chi(n), ~~~ q_n =
1-\xi_n -\chi(n) Y_n, ~~~ n \geq n_0.
\eea
Note that our conditions
on the variables in (\ref{1006b})
ensure that $(\eps/2) \leq p_n \leq 1 -(\eps/2)$
almost surely for all $n$, so that for a.e.~$\omega$,
$p_n$ and $q_n$ are true probabilities
bounded strictly away from $0$ and from $1$.

For $n=1,2,\ldots$, we
set
\bea
\label{0520b}
\zeta_n := \log \left( \frac{\xi_n}{1-\xi_n} \right).
\eea
Write $\bbE$ for expectation under $\bbP$.

In our model, by (\ref{chicon}), $\chi(n) \to 0$ as $n \to
\infty$. Thus, from (\ref{1006c}), in the limit $n \to \infty$, we
approach the well-known random walk in i.i.d.~random environment
as studied in \cite{kozlov}, \cite{solomon} and subsequently. In
addition, when $\bbE[\zeta_1]=0$, in the limit as $n \to \infty$
we approach the critical case often referred to as {\em Sinai's
regime} after \cite{sinai}. Our results show that despite this,
the behaviour of our model is, in general, very different to the
behaviour of these limiting cases, depending on the nature of the
perturbation $\chi$.

In work in preparation, we study the long-run limiting
behaviour (as $t \to \infty$) of our random walk $\eta_t(\omega)$
in terms of its distance from the origin. Of interest are both the
almost sure and `in probability' (see, for example,
\cite{sinai,cp}) behaviour. In Sinai's regime for the RWRE on
$\Z^+$, Comets, Menshikov and Popov (\cite{cmp}, Theorem 3.2) show
that, for a.e.~$\omega$ and any $\eps>0$,
\[ \frac{\eta_t(\omega)}{(\log t)^2} < (\log \log t)^{2+\eps}, ~{\rm a.s.}\]
for all but finitely many $t$ (where a.s.~stands for $\Pr$-almost surely). 
This result (for the RWRE on $\Z$)
dates back to Deheuvels and Rev\'esz \cite{dere}. An exact upper
limit result is given in \cite{hushi}. In work in preparation, we
study analogous almost sure results (in both null-recurrent and
transient cases) for our perturbed RWRE. For example, in the
$\bbP$-almost sure transient case of the RWRE perturbed from
Sinai's regime (that is, with $\chi(n)=n^{-\alpha}$ for some fixed $0<\alpha
<1/2$, we have $\bbE[\zeta_1]=0$, $\Var[\zeta_1]>0$ and
$\lambda < 0$, where $\lambda$ is defined at (\ref{0530a})), 
we have that for a.e.~$\omega$, for any $\eps>0$, as $t
\to \infty$, \bean
 (\log \log \log t)^{-(1/\alpha)-\eps} <
\frac{\eta_t(\omega)}{(\log t)^{1/\alpha}} < (\log \log
t)^{(2/\alpha)+\eps} , ~{\rm a.s.} \eean for all but finitely many $t$. Thus in
this case, we see that the random walk, for almost every
environment, is contained in a window about $(\log t)^{1/\alpha}$.
This aspect of the problem requires additional techniques,
however, and we do not discuss this further in the present paper.

In the next section we state our results. Theorems \ref{iid}, \ref{0928t}, and \ref{thm1022}
are special cases of the model in which some of the
random variables $\xi_i$ and $Y_i$ are degenerate (that is, equal to a constant almost surely). In particular,
Theorems \ref{iid} and \ref{0928t} include some known results, when our model reduces
to previously studied systems. In Theorem \ref{thm2}, the underlying environment
is not in the `critical regime'. Our main results, Theorems
\ref{cor2} and \ref{cor3},
deal with the main case of interest, in which the underlying environment is truly random and is,
in a sense to be demonstrated, critical.

\section{Main results}

Most of our results
will be formulated for almost all environments $\omega$ (in some
sense, for all `typical' environments), that is $\bbP$-almost surely over
$(\Omega, \F, \bbP)$.

If $Y_1 =0$ $\bbP$-a.s., then our model reduces to
the standard reflected one-dimensional
random walk in an i.i.d.~random environment. In this case $p_n = \xi_n$
and $q_n=1-\xi_n$, $n=1,2,\ldots$, and so (with the definition
at (\ref{0520b})) $\zeta_n= \log(p_n/q_n)$.
Criteria for recurrence of the
 RWRE
$\eta_t(\omega)$ in this case
were given by Solomon
\cite{solomon}, for the case
where $( \xi_i ; i=1,2,\ldots)$ is an i.i.d.~sequence, and generalised by Alili \cite{alili}.
For the case in which larger
jumps are permitted, see, for example, \cite{key}.

The following well-known
result dates back to
Solomon \cite{solomon}.
\begin{theorem}
\label{iid}
Let $(\eta_t(\omega); t \in \Z^+)$
be the random walk in i.i.d.~random
environment, with $\bbP[Y_1=0]=1$.
Suppose $\Var [ \zeta_1]>0$.
\begin{itemize}
\item[(i)] If $\bbE [ \zeta_1 ] <0$, then $\eta_t(\omega)$ is
transient for a.e.~$\omega$. \item[(ii)] If $\bbE [ \zeta_1 ] = 0$, then
$\eta_t(\omega)$ is null-recurrent for a.e.~$\omega$. \item[(iii)] If
$\bbE [ \zeta_1 ] >0$, then $\eta_t(\omega)$ is
ergodic for a.e.~$\omega$.
\end{itemize}
\end{theorem}
The critical (null-recurrent) regime $\bbE[ \log(p_1/q_1) ]=0$ is
known as {\em Sinai's regime}, after \cite{sinai}. This regime has
been extensively studied; see, for example,
\cite{cp,hushi,ktt,kks}. For an outline proof of Theorem \ref{iid}
using Lyapunov function methods, similar to those employed in this
paper, see Theorem 3.1 of \cite{cmp}. In this paper we extend the
classification criteria of Theorem \ref{iid} to encompass the case
in which the $p_n$ are not i.i.d.~and in which $\bbE[
\log(p_n/q_n) ]$ is {\em asymptotically} zero, as $n \to \infty$.
Our results are, in some sense, a random environment analogue of
those for Markov processes with asymptotically zero mean drift
given in \cite{mai} (see below).

For the remainder of the paper we suppose $\bbP[Y_1=0] <1$. This includes the interesting case
where $Y_1 = b$ $\bbP$-a.s., for some $b \in [-1,1] \setminus \{ 0 \}$.
Our techniques do, however, enable us
to allow $Y_1$ to be random.

Although not as famous as the RWRE, another
system that has been well studied is the rather classical
problem of a Markov chain with asymptotically zero drift. This problem
was studied by Lamperti \cite{lamp1,lamp2}. General criteria for recurrence,
transience and ergodicity were given by Menshikov, Asymont, and Iasnogorodskii in \cite{mai}.
Theorem \ref{0928t} below is a consequence of their main result,
Theorem 3, applied to our problem when $\Var[ \zeta_1]=0$ and
$\Var[Y_1]=0$; that is, the distributions of $\xi_1$ and $Y_1$ are
both degenerate (i.e.~equal to a constant almost surely). In
particular, we have a {\em non-random} environment $\omega$. If,
on the other hand,
 $\xi_1$ is degenerate but $Y_1$ is not,
then we have a random (asymptotically small) perturbation on an
underlying non-random environment, and we have Theorem
\ref{thm1022} below.

We use the notation $\log_1 x := \log x$ and
$\log_k x := \log (\log_{k-1} x)$ for $k=2,3,\ldots$.
\begin{theorem}
\label{0928t}
Suppose $\bbP[Y_1=b]=1$ for some $b \in [-1,0) \cup (0,1]$.
Suppose $\bbP [ \xi_1 = c]=1$ for some $c \in (0,1)$.
\begin{itemize}
\item[(i)] If $c<1/2$, then $\eta_t(\omega)$ is transient.
\item[(ii)] If $c>1/2$, then $\eta_t(\omega)$ is ergodic.
\item[(iii)] Suppose $c=1/2$. Suppose there exist $s \in \Z^+$ and
$K \in \N$ such that, for all $n \in [K, \infty)$ and some $h>1$
the following inequality holds: \bea \label{0928ta} b \chi(n) >
\frac{1}{4n} + \frac{1}{4n \log n} + \cdots + \frac{h}{4n
\prod_{i=1}^s \log_i n}. \eea Then $\eta_t(\omega)$ is ergodic.
\item[(iv)] Suppose $c=1/2$. Suppose there exist $s,t \in \Z^+$
and $K \in \N$ such that, for all $n \in [K, \infty)$ and some
$h<1$ the following inequality holds: \bea \label{0928tb}
-\frac{1}{4n} - \frac{1}{4n \log n} - \cdots - \frac{h}{4n
\prod_{i=1}^s \log_i n} \leq b \chi(n) \nonumber\\
\leq
\frac{1}{4n} + \frac{1}{4n \log n} + \cdots + \frac{h}{4n
\prod_{i=1}^t \log_i n} . \eea Then $\eta_t(\omega)$ is
null-recurrent. \item[(v)] Suppose $c=1/2$. Suppose there exist $s
\in \Z^+$ and $K \in \N$ such that, for all $n \in [K, \infty)$
and some $h>1$ the following inequality holds: \bea \label{0928tc}
b \chi(n) \leq -\frac{1}{4n} - \frac{1}{4n \log n} - \cdots -
\frac{h}{4n \prod_{i=1}^s \log_i n}. \eea Then $\eta_t(\omega)$ is
transient.
\end{itemize}
\end{theorem}
Theorem \ref{0928t} follows directly by applying Theorem 3 of
\cite{mai} to our case, with $m(x)=-2\chi(x)$ and $b(x)=1$.
\\

\remark In the case $c=1/2$ the critical case in terms of the
recurrence, transience and ergodicity is when the perturbation
$\chi(n)$ is, ignoring logarithmic terms, of order $n^{-1}$; we
say that the `critical exponent' is $-1$. This contrasts with our
results in the case where $\Var[\xi_1]>0$ (see Theorems \ref{cor2}
and \ref{cor3}),
in which the critical exponent is $-1/2$. \\

The following result deals with the case in which the distribution
of $\xi_1$ is degenerate, but that of $Y_1$ is not; in this case
we have a homogeneous non-random environment subject to an
asymptotically small random perturbation. In particular, parts
(iii) and (iv) of the theorem deal with the case when the
underlying environment is that of the simple random walk. Here,
$\eqd$ stands for equality in distribution.
\begin{theorem}
\label{thm1022} Suppose $\bbP [ \xi_1 = c]=1$ for some $c \in
(0,1)$, and $\Var[Y_1]>0$.
\begin{itemize}
\item[(i)] If $c<1/2$, then $\eta_t(\omega)$ is
transient for a.e.~$\omega$. \item[(ii)] If $c>1/2$, then
$\eta_t(\omega)$ is ergodic for a.e.~$\omega$. \item[(iii)] If $c=1/2$
and $Y_1 \eqd -Y_1$, then $\eta_t(\omega)$ is
null-recurrent for a.e.~$\omega$. \item[(iv)] Suppose $c=1/2$ and
$\bbE[Y_1] \neq 0$. Suppose $\chi(n)=an^{-\beta}$ for $a>0$ and
$\beta>0$.
\begin{itemize}
\item[(a)] If $0< \beta <1$ and $\bbE[Y_1] >0$ then
$\eta_t(\omega)$ is ergodic for a.e.~$\omega$. \item[(b)] If $\beta >1$
then $\eta_t(\omega)$ is null-recurrent for a.e.~$\omega$. \item[(c)] If
$0<\beta <1$ and $\bbE[Y_1] <0$ then $\eta_t(\omega)$ is
transient for a.e.~$\omega$.
\end{itemize}
\end{itemize}
\end{theorem}
We prove Theorem \ref{thm1022} along with our main results in
Section \ref{prfs}. \\

\remarks Note that in part (iii), $Y_1 \eqd -Y_1$ implies that all
odd moments of $Y_1$ are zero. By modifications to the proof
of Theorem \ref{thm1022} one can obtain a more refined result,
specifically that with $p:= \min \{ j \in \{1,3,5,\ldots\} :
\bbE[Y_1^j] \neq 0\}$, for $p>1$ we have a statement analogous to
part (iv) but with $\bbE[Y_1]$ replaced by $\bbE[Y_1^p]$ and with
the critical value of $\beta$ being $1/(2(p-1))$ for $p>1$, rather
than $1$. We do not go into further detail here.

 Theorem \ref{thm1022} (iv) demonstrates that in the case of a randomly perturbed
simple random walk, the critical exponent for the perturbation
is $-1$, as in the case of the non-random perturbation (Theorem \ref{0928t}). It may be
possible to refine Theorem \ref{thm1022} (iv) to obtain more delicate results analogous
to those of Theorem \ref{0928t}. \\

For the remainder of the paper,
we ensure that the underlying
environment is {\em random}, by
supposing $\Var [ \zeta_1] >0$.
First we consider the case $\bbE[\zeta_1] \neq 0$. Here we have Theorem \ref{thm2} below.
In this situation, the perturbation introduced by $\chi(n)Y_n$
does not affect the criteria given in (i) and (iii) of Theorem \ref{iid}.
\begin{theorem}
\label{thm2} Suppose $\Var[ \zeta_1]>0$, $\bbE[\zeta_1] \neq 0$,
and $\bbP[Y_1=0] <1$.
\begin{itemize}
\item[(i)] If $\bbE [ \zeta_1 ] <0$, then  $\eta_t(\omega)$ is
transient for a.e.~$\omega$. \item[(ii)] If $\bbE [ \zeta_1 ] >0$, then
$\eta_t(\omega)$ is ergodic for a.e.~$\omega$.
\end{itemize}
\end{theorem}
The proof of the theorem follows using the same methods as employed
in the proof of Theorem 3.1 of \cite{cmp} or later in this paper, but is 
essentially
simpler than for our main results. 
We can construct a `martingale' (as at (\ref{0818b}) below) which is easily shown 
(by the Law of the Iterated Logarithm, Lemma \ref{itlog} below)
to be bounded or tend to infinity for a.e.~$\omega$. Similarly for the stationary
measure. The theorem then follows by our Lyapunov function criteria (Lemmas \ref{crit1}
and \ref{crit2} below). We follow this method in detail, in less straightforward cases, later in the paper,
and so do not repeat the argument here.

For the remainder of the paper we consider the more interesting
case where $\bbE[\zeta_1]=0$, so that we have a random walk in a
random environment that is asymptotic to Sinai's regime. We prove
general results about this RWRE with asymptotically zero
perturbation that are analogous to Theorem \ref{0928t}, but
significantly different.

If $\bbP[Y_1=0]<1$ (and
permitting the case that $\bbP[Y_1=c]=1$ for
some $c$ with $0<|c| \leq 1$) we define
\bea
\label{0530a}
 \lambda := \bbE \left[ \frac{Y_1}{\xi_1 (1-\xi_1)} \right] .
\eea Also, we use the notation \bea \label{1112a} \sigma^2 :=
\Var[ \zeta_1 ]. \eea Note that, under the condition (\ref{ue}),
we have $\sigma^2<\infty$ and, since $Y_1$ is bounded,
$|\lambda|<\infty$. We also draw attention to the fact that, given
(\ref{ue}), $\bbP$-a.s., \bea \label{bndd} -\frac{1}{\eps^2} \leq
\frac{Y_1}{\xi_1 (1-\xi_1)} \leq \frac{1}{\eps^2}, \eea
a fact that we shall use later. For what follows, of separate
interest are the two cases $\lambda =0$ and $\lambda \neq 0$. We
concentrate on the latter case for most of the results that follow
(but see the remark after Theorem \ref{cor3}). However, our first
result deals with the case in which $Y_1/\xi_1 \eqd
-Y_1/(1-\xi_1)$. This implies $\lambda=0$ (see (\ref{0530a})), but
is a rather special case; Theorem \ref{cor0} demonstrates that in
this case the detailed behaviour of $\chi$ is not important: as
long as $\chi(n) \to 0$ as $n \to \infty$, then $\eta_t (\omega)$
is null-recurrent for a.e.~$\omega$.
\begin{theorem}
\label{cor0} With $\sigma$ as defined at (\ref{1112a}), suppose
that $Y_1/\xi_1 \eqd -Y_1/(1-\xi_1)$, $\bbP[Y_1=0]<1$,
$\bbE[\zeta_1]=0$, and $\sigma^2>0$. Then $\eta_t(\omega)$ is
null-recurrent for a.e.~$\omega$.
\end{theorem}
An example of $(Y_1,\xi_1)$ for which Theorem \ref{cor0} holds
is when $Y_1$ and $\xi_1$ are independent uniform
random variables on $(-1,1)$ and $(\eps, 1-\eps)$ respectively.

Our remaining results deal with the case $\lambda \neq 0$ (but see
also the remark after Theorem \ref{cor3}). In our next result
(Theorem \ref{cor2}), we give some rather specific conditions on
the asymptotic behaviour of the function $\chi$. Theorem
\ref{cor2} is a special case of our general result, Theorem
\ref{cor3}.
\begin{theorem}
\label{cor2} With $\lambda$ and $\sigma$ defined at (\ref{0530a})
and (\ref{1112a}) respectively, suppose that $\lambda \neq 0$,
$\bbP[Y_1=0]<1$, $\bbE[\zeta_1]=0$, and $\sigma^2>0$. Let $c_{\rm
crit}:=\sigma 2^{-1/2}$.
\begin{itemize}
\item[(i)] If there exist constants $c>c_{\rm crit}$ and $n_0 \in
\Z^+$ such that $\lambda \chi(n) \geq cn^{-1/2}(\log \log
n)^{1/2}$ for all $n \geq n_0$, then $\eta_t(\omega)$ is
ergodic for a.e.~$\omega$. 
\item[(ii)] If there exist constants $c \leq
c_{\rm crit}$ and $n_0 \in \Z^+$ such that $|\lambda| \chi(n) \leq
cn^{-1/2}(\log \log n)^{1/2}$ for all $n \geq n_0$, then
$\eta_t(\omega)$ is null-recurrent for a.e.~$\omega$. 
\item[(iii)] If
there exist constants $c>c_{\rm crit}$ and $n_0 \in \Z^+$ such
that $\lambda \chi(n) \leq -cn^{-1/2}(\log \log n)^{1/2}$ for all
$n \geq n_0$, then $\eta_t(\omega)$
 is transient for a.e.~$\omega$.
\end{itemize}
\end{theorem}
\remark
\label{rem3}
 Theorem \ref{cor2} shows that in our case the critical
exponent for the perturbation is $-1/2$. This contrasts with the
deterministic environment case (as in Theorem \ref{0928t}, and see
\cite{mai}, Theorem 3), in which the critical exponent is $-1$.
When the perturbation is smaller than this critical size (as in
part (ii)), it is insufficient to change the recurrence/transience
characteristics of the Markov chain from those of Sinai's regime.
If the perturbation is greater
 than the critical size, it changes the behaviour
 of the
Markov chain from that of Sinai's regime, making it either transient or ergodic depending
on the sign of the perturbation. This feature is present in our most general result, Theorem \ref{cor3}. \\

Theorem \ref{cor2} will follow as a corollary to Theorem
\ref{cor3}, below. Theorem \ref{cor3} is more refined than Theorem
\ref{cor2}.
In order to formulate our deeper result,
we need more precise conditions on the behaviour of the
perturbation function $\chi(n)$. To achieve this, we define the
notions of $k$-{\em supercritical} and $k$-{\em subcritical}
below. First, we need some additional notation.

Recall the notation $\log_1 (x) := \log (x)$ and $\log_k(x) :=
\log (\log_{k-1} (x))$ for $k=2,3,\ldots$. Let $n_k$ denote the
smallest positive integer such that $\log_{k+1}(n_k) \geq 0$. Let
$a_k := 2$ for $ k \in \N \setminus \{ 3\}$ and $a_3 :=3$. For
each $k \in \N$ we define the $[0,\infty)$-valued function
$\varphi_k$ as follows (we use the given form for the $\varphi_k$
due to the appearance in the sequel of the Law of the Iterated
Logarithm). For $x \in [e,\infty)$ and $d \in \R$, let
\[ \varphi_1(x;d) := ((2+d)\log_2{x})^{1/2} ,\]
and for $k=2,3,\ldots$, with
$x \in [n_k,\infty)$
and $d \in \R$, let
\bea
\label{0604a}
\varphi_k(x;d) := \left( \sum_{i=1}^{k-1} a_{i+1} \log_{i+1} x
+ ( a_{k+1} +d ) \log_{k+1} x \right)^{1/2}.
\eea
We shall see that the behaviour of the Markov chain
$\eta_t(\omega)$ is determined by the driving function $\chi$. By
applying the Law of the Iterated Logarithm, we shall see that the
critical form of $\chi$ is related to an iterated logarithm
expression of the form of $\varphi_k$.

In order to formulate our main result
we make the following definitions of $k$-{\em supercritical} and
$k$-{\em subcritical}.
\begin{definition}
\label{def1}
Recall the definitions of $\lambda$ and $\sigma$
at (\ref{0530a}) and (\ref{1112a}) respectively.
Suppose $\lambda \neq 0$. For $k \in\N$, we say $\chi$ is $k$-{\em supercritical} if
there exist constants $c \in (0,\infty)$ and
$n_0 \in \Z^+$, such that, for all $n \geq n_0$,
\bea
\label{super}
\chi(n) \geq \frac{\sigma}{ 2|\lambda|} n^{-1/2} \varphi_k (n;c).
\eea
For $k \in \N$, we say
$\chi$ is $k$-{\em subcritical} if
there exist constants $c \in (0,\infty)$ and
$n_0 \in \Z^+$ such that, for all $n \geq n_0$,
\bea
\label{sub}
\chi(n) \leq \frac{\sigma}{ 2|\lambda|} n^{-1/2} \varphi_k (n;-c).
\eea
\end{definition}
\remarks Implicit in $\chi$ being $k$-subcritical or
$k$-supercritical is the constant $c$, a fact that we make
repeated use of in the proofs in Section \ref{prfs}. Whenever we
consider a $k$-subcritical or $k$-supercritical function in what
follows, we understand this to imply the existence of such a $c$,
and often refer to the constant $c$ in this context.

Also, observe that if for some $k \in \N$, $\chi$ is $k$-supercritical,
with implicit constant $c \in (0,\infty)$, then for any $c' \in (0,c)$ we have that
(\ref{super}) implies
\[ \chi(n) \geq \frac{\sigma}{ 2|\lambda|} n^{-1/2} \varphi_k (n;c) \geq
\frac{\sigma}{ 2|\lambda|} n^{-1/2} \varphi_k (n;c').\]
Similarly if for some $k \in \N$, $\chi$ is $k$-subcritical,
with implicit constant $c \in (0,\infty)$, then for any $c' \in (0,c)$ we have that
(\ref{sub}) implies
\[ \chi(n) \leq \frac{\sigma}{ 2|\lambda|} n^{-1/2} \varphi_k (n;-c) \leq
\frac{\sigma}{ 2|\lambda|} n^{-1/2} \varphi_k (n;-c').\]

Finally, we note that Definition \ref{def1} excludes functions
that oscillate
significantly about the critical region $n^{-1/2}$. \\

Our most general result is as follows.
\begin{theorem}
\label{cor3}  With $\lambda$ and $\sigma$ defined at (\ref{0530a})
and (\ref{1112a}) respectively, suppose that $\lambda \neq 0$,
$\bbP[Y_1=0]<1$, $\bbE[\zeta_1]=0$ and $\sigma^2>0$.
\begin{itemize}
\item[(i)] If, for some $k\in\N$, $\chi$ is $k$-supercritical
(\ref{super}) and $\lambda >0$, then $\eta_t(\omega)$
 is ergodic for a.e.~$\omega$.
\item[(ii)]
If, for some $k\in\N$, $\chi$ is $k$-subcritical (\ref{sub})
 then $\eta_t(\omega)$
is null-recurrent for a.e.~$\omega$. 
\item[(iii)] If, for some $k\in\N$,
$\chi$ is $k$-supercritical (\ref{super}) and $\lambda <0$, then
$\eta_t(\omega)$
 is transient for a.e.~$\omega$.
\end{itemize}
\end{theorem}
\remark In the general case $\lambda=0$, it turns out that higher
moments contribute, and we obtain a slightly more general form of
Theorem \ref{cor3}. It is straightforward to modify the proof of
Theorem \ref{cor3} to obtain such a result. Specifically, if for
$r \in \N$ we set
\[ \lambda_r := \frac{1}{r} \bbE \left[ Y_1^r \left( \frac{1}{(1-\xi_1)^r} + \frac{(-1)^{r+1}}{\xi_1^r}
\right) \right],\] and $p:=\min \{ j \in \N : \lambda_j \neq 0\}$,
then for $p>1$ a statement of the form of Theorem \ref{cor3} holds
but with $\lambda$ replaced by $\lambda_p$ and the conditions on
$\chi$
being replaced by conditions on $\chi^p$. We do not pursue the details here. \\

We will prove Theorem \ref{cor3} in the next section. The idea
behind the proof of the recurrence and transience conditions is to
construct a function $f$ of the process $\eta_t(\omega)$ such that
$f(\eta_t(\omega))$ is a `martingale' everywhere except in a
finite region, and determine the cases in which this function is
finite or infinite. The proof of ergodicity relies on the
construction of a stationary measure and determining its
properties.

\section{Proofs of main results}
\label{prfs}
Before embarking upon the proof of Theorem \ref{cor3},
we need some preliminary results.
First we present the criteria for classification of countable
Markov chains which we will require.

Let $(W_t; t \in \Z^+)$ be a discrete, irreducible, aperiodic, time-homogeneous Markov chain
on $\Z^+$.
We have the following
classification criteria, which are
consequences of those given in Chapter 2 of \cite{fmm}. The following
result, which we state without proof,
 is a consequence of Theorem 2.2.2 of \cite{fmm}, and is
slightly more general than Proposition 2.1 of \cite{cmp}.
\begin{lemma}
\label{crit1}
Suppose there exist a function $f: \Z^+ \to [ 0 , \infty)$ which
is uniformly bounded and nonconstant, and a set $A \subset \Z^+$
such that
\bea
\label{0928a}
\Exp [ f( W_{t+1}) - f(W_t) | W_t = x ] =0,
\eea
for all $x \in \Z^+ \setminus A$, and
\bea
\label{0928b}
f(x) > \sup_{y \in A} f(y),\eea
for at least one $x \in \Z^+ \setminus A$. Then
the Markov chain $(W_t)$ is transient.
\end{lemma}
The following result is contained in Theorem 2.2.1 in \cite{fmm}.
\begin{lemma}
\label{crit2}
 Suppose that there exist a function $f: \Z^+ \to [0,\infty)$
and a finite set $A \subset \Z^+$
such
that
\bea
\label{cr2}
 \Exp \left[ f( W_{t+1})- f(W_t) | W_t = x \right] \leq 0,
\eea
for all $x \in \Z^+ \setminus A$, and
$f(x) \to +\infty$ as $x \to \infty$.
Then the Markov chain $(W_t)$ is recurrent.
\end{lemma}
We will need Feller's refined form for the Law of the Iterated
Logarithm \cite{feller}. The following result is a consequence of
Theorem 7 of \cite{feller}.
\begin{lemma}
\label{itlog} Let $X_i$, $i=1,2,\ldots$, be a sequence of
independent random variables with $\Exp[X_i]=0$ for all $i$, and
$\Exp[X_i^2]=\sigma_i^2 < \infty$ for $i=1,2,\ldots$. Suppose the
$X_i$ are bounded, that is $\Pr[|X_i| > C] =0$ for all $i$ and
some $0<C<\infty$. Let \bea \label{1005a} s_n^2 := \sum_{i=1}^n
\sigma_i^2 .\eea Suppose that $s_n \to \infty$ as $n \to \infty$.
  Let $S_n := \sum_{i=1}^n X_i$.
For some $k \in \N$ and $\eps \in (-\infty, \infty)$,
define $\varphi_k(n;\eps)$ as at (\ref{0604a}).
Then
\bea
\label{0603b}
 \Pr \left[ S_n > s_n \varphi_k( s_n^2;\eps) ~\rm{i.o.} \right]
= \left\{ \begin{array}{ll} 1 & {\rm if}~ \eps<0 \\
0 & {\rm if}~ \eps >0 \end{array} \right. \eea In particular, if
the $X_i$ are i.i.d.~and bounded random variables with
$\Exp[X_1^2]=\sigma^2$,
we have \bea \label{0603bb}
 \Pr \left[ S_n > \sigma n^{1/2} \varphi_k(n;\eps) ~\rm{i.o.} \right]
= \left\{ \begin{array}{ll} 1 & {\rm if}~ \eps<0 \\
0 & {\rm if}~ \eps >0 \end{array} \right.
\eea
\end{lemma}
We will also need the following result. Recall the definition of
$\varphi_k(i;d)$ at (\ref{0604a}).
\begin{lemma}
\label{0604c} For $k \in \N$, let $n_k$ be the smallest positive
integer such that $\log_{k+1} n_k \geq 0$. For any $d \in \R$, we
have \bea \label{0604b} \sum_{i=n_k}^{n} i^{-1/2} \varphi_k(i;d) =
2n^{1/2} \varphi_k(n;d) + \alpha_n, \eea where $|\alpha_n| <
6n^{1/2}$ for all $n$ sufficiently large.
\end{lemma}
\proof We have, for $k \in \N$, \[ \frac{\ud}{\ud x} \left(
x^{1/2} \varphi_k(x;d) \right) = \frac{1}{2} x^{-1/2}
\varphi_k(x;d) + x^{1/2} \varphi_k'(x;d) , \] where \[
\varphi_k'(x;d) = \frac{1}{2} \left(\varphi_k (x ; d) \right)^{-1}
\cdot \left( \frac{2}{x \log x} + \frac{3}{x \log {x} \log
\log{x}} + \cdots \right) < \frac{1}{x} ,\] for $x$ sufficiently
large. Thus, for any $k \in \N$, \bea \label{0607a} \int_{n_k}^n
x^{-1/2} \varphi_k(x;d) \ud x & = & 2\left[ x^{1/2} \varphi_k(x;d)
\right]_{n_k}^n
-2 \int_{n_k}^n x^{1/2} \varphi_k' (x;d) \ud x \nonumber\\
& = & 2 n^{1/2} \varphi_k(n;d) + b_n,
\eea
where
\[ |b_n| \leq  2 \int_{n_k}^n x^{1/2} \varphi_k' (x;d) \ud x +2n_k^{1/2} \varphi_k
(n_k;d)
 \leq C_k + 2 \int_0^n x^{-1/2} \ud x,\]
for some $0<C_k<\infty$, which depends on $k$ (and $d$). Thus, for
each $k$, $|b_n| \leq 5n^{1/2}$ for all $n$ sufficiently large.
Since $x^{-1/2} \varphi_k(x;d)$ is a decreasing function for all
$x$ sufficiently large (depending on $k$ but not $d$), we have
that there exist finite positive constants $C_k'$ and $C''_k$ such
that \bean \sum_{i=n_k}^{n} i^{-1/2} \varphi_k(i;d) +C'_k \geq
\int_{n_k}^{n} x^{-1/2} \varphi_k(x;d) \ud x \geq
\sum_{i=n_k+1}^{n} i^{-1/2} \varphi_k(i;d)-C''_k. \eean So we have
\bea \label{0921b} 0 \leq \sum_{i=n_k}^{n} i^{-1/2} \varphi_k(i;d)
- \int_{n_k}^{n} x^{-1/2} \varphi_k(x;d) \ud x \leq n_k^{-1/2}
\varphi_k(n_k;d) +C, \eea for some $0<C<\infty$, that does not
depend on $n$. Then from (\ref{0921b}) and (\ref{0607a}) we obtain
(\ref{0604b}). $\square$ \\

For a given realization $\omega$ of our random environment, with
$p_i$ and $q_i$, $i=1,2,\ldots$ defined by (\ref{1006b}), let \bea
\label{0520a} D(\omega) := \sum_{i=1}^\infty \frac{1}{q_i}
\prod_{j=1}^i \frac{q_j}{p_j} = \frac{1}{p_1} + \frac{q_1}{p_1
p_2} + \frac{q_1 q_2}{p_1 p_2 p_3} + \cdots . \eea
\begin{lemma}
\label{stat} If, for a given environment $\omega$, the quantity
$D(\omega)$ as defined at (\ref{0520a}) is finite, then the Markov
chain $\eta_t(\omega)$ is ergodic. On the other hand, if
$D(\omega)=\infty$, then the Markov chain $\eta_t(\omega)$ for
this $\omega$ is not ergodic.
\end{lemma}
\proof For fixed environment $\omega$, i.e., given a configuration
of $(p_i ; i=1,2,\ldots)$, $\eta_t(\omega)$ is a reversible Markov
chain. For this Markov chain one has the stationary measure $\mu =
(\mu_0, \mu_1, \ldots )$, where \[
\mu_0 = 2, ~
\mu_1 = \frac{1}{p_1}, ~ {\rm and} ~ \mu_n = \frac{1}{p_1}
\prod_{i=1}^{n-1} \frac{q_i}{p_{i+1}}, ~~~ n \geq 2. \] Then,
with the definition of $D(\omega)$ at (\ref{0520a}), we have
\[ \sum_{i=0}^\infty \mu_i = 2 + D(\omega).\]
Thus, if, for this $\omega$, $D(\omega)$ is finite, then the
Markov chain $\eta_t(\omega)$ is ergodic, since we can obtain a
stationary distribution. On the other hand, if $D(\omega)=\infty$
for this $\omega$, the Markov chain $\eta_t(\omega)$
 is not ergodic. $\square$ \\

Our next result,
Lemma \ref{0530b}, uses the Law of the Iterated
Logarithm to analyse the
behaviour of sums of i.i.d.~random variables
weighted by the function $\chi$.
\begin{lemma} \label{0530b}
Let $Z_i$, $i=1,2,\ldots$, be a sequence of i.i.d.~random
variables which are bounded (so that $\Pr[|Z_1|>B]=0$ for some
$0<B<\infty$), such that $\Exp[Z_1] \geq 0$. Let $\chi:
[0,\infty) \to [0,\infty)$ such that (\ref{chicon}) holds. 
With $\lambda$ defined at (\ref{0530a}), suppose $\lambda
\neq 0$.
\begin{itemize}
\item[(a)]
Suppose $\Exp[Z_1]>0$.
 Suppose that, for some $k\in\N$, $\chi$ is $k$-subcritical as defined at (\ref{sub}).
Then with probability one, for any $\eps>0$, for all but finitely
many $n$, \bea \label{0922r} -n^{\eps} \leq \sum_{i=1}^n Z_i
\chi(i) \leq \frac{\sigma \Exp[Z_1]}{|\lambda|} n^{1/2} \varphi_k
(n; -c/3). \eea \item[(b)] Suppose $\Exp[Z_1]>0$. Suppose that,
for some $k\in\N$, $\chi$ is $k$-supercritical as defined at
(\ref{super}). Then with probability one, for all but finitely
many $n$, \bea \label{0922t} \sum_{i=1}^n Z_i \chi(i) \geq
\frac{\sigma \Exp[Z_1]}{|\lambda|} n^{1/2} \varphi_k ( n; c/3).
\eea \item[(c)] Suppose $\Exp[Z_1]=0$. Then for any $\eps>0$ with
probability one, for all but finitely many $n$, \bea \label{1007a}
\sum_{i=1}^n Z_i \chi(i) \leq \eps ( n \log \log n)^{1/2}. \eea
\end{itemize}
\end{lemma}
\remark When we come to apply Lemma \ref{0530b} later in the proofs of the
theorems, the configuration $(Z_i, i \geq 1)$ that we will use
will be specified by the realization of the random
environment $\omega$, so that the qualifier `with probability one'
in the lemma translates as `for a.e.~$\omega$.' \\

\noindent {\bf Proof of Lemma \ref{0530b}.} Recall the definitions
of $\lambda$ and $\sigma$ at (\ref{0530a}) and (\ref{1112a})
respectively. Suppose $\lambda \neq 0$. For the proofs of parts
(a) and (b), suppose that $\Exp[Z_1]>0$. First we prove part (a).
Suppose that for some $k \in \N$ $\chi$ is $k$-subcritical. Write
\bea \label{1022a} S_n := \sum_{i=1}^n (Z_i-\Exp[Z_i]) \chi(i) .
\eea Then \bea \label{1216a}
 \Var [S_n] = \Var[Z_1] \sum_{i=1}^n (\chi(i))^2.\eea
Suppose that $\Var[S_n] \to \infty$ as $n \to \infty$.
Then, by Lemma \ref{itlog}, taking
$X_i=(Z_i-\Exp[Z_i])$, we have that
with probability one the configuration of $(Z_i, i \geq 1)$ is such that
\[ S_n > ( \Var[S_n])^{1/2} (3 \log \log{ (\Var[S_n])} )^{1/2} ,\]
for only finitely many $n$. (The constant $3$ appears for the sake of simplicity,
any constant strictly greater than $2$ will suffice). That is, with probability one,
for all but finitely many $n$,
\[ S_n \leq ( \Var[S_n])^{1/2} (3 \log \log{ (\Var[S_n])} )^{1/2}
\leq ( \Var[S_n])^{1/2} (3 \log \log n )^{1/2},\] the second
inequality following from (\ref{1216a}) and (\ref{sub}). Thus,
with probability one, using (\ref{1216a}) and (\ref{sub}) once
more, we have that for any $\eps>0$ and for all but finitely many
$n$, $S_n \leq n^\eps$. Thus, with probability one, for all but
finitely many $n$, since $\Exp[Z_1]>0$ and $\chi$ is a nonnegative
function, \bea \label{0922aa} -n^{\eps} \leq \sum_{i=1}^n Z_i
\chi(i) \leq n^\eps + \Exp[Z_1] \sum_{i=1}^n \chi(i). \eea The
lower bound in (\ref{0922aa}) establishes the lower bound in
(\ref{0922r}). We now need to prove the upper bound. By
(\ref{sub}), we have that there exist $c \in (0,\infty)$ and $k
\in \N$ such that for all $n$ sufficiently large \bea
\label{0922aab} \sum_{i=1}^n \chi(i)  \leq  \frac{\sigma}{2
|\lambda|} \sum_{i=1}^n i^{-1/2} \varphi_k (i; -c/2) .\eea Then
from (\ref{0922aab}) with (\ref{0604b}) we obtain, for all $n$
sufficiently large \bea \label{0922aabx} \sum_{i=1}^n \chi(i)
 \leq \frac{\sigma}{|\lambda|} n^{1/2} \varphi_k(n;-c/2) + \frac{3 \sigma}{|\lambda|} n^{1/2}.
\eea Hence from (\ref{0922aabx}) and the upper bound in
(\ref{0922aa}), we have that, with probability one, for all but
finitely many $n$,
\[
\sum_{i=1}^n Z_i \chi(i) \leq \frac{ \sigma \Exp[Z_1]}{|\lambda|}
n^{1/2} \varphi_k (n; -c/2) + \frac{3 \sigma \Exp[Z_1]}{|\lambda|}
n^{1/2}+n^{\eps} .\] Then we can absorb the final two terms on the
right hand side to give (\ref{0922r}), given that $\Var[S_n] \to
\infty$ as $n \to \infty$. On the other hand, suppose that
$\Var[S_n] \leq C$ for all $n$ and some $C < \infty$. Then, by
(\ref{1216a}), we have that $\sum_{i=1}^n (\chi(i))^2 < C$ for
some $0<C<\infty$. So, by Jensen's inequality, and the boundedness
of the $Z_i$, we have that for all $n$, \bea \label{987a}
\sum_{i=1}^n Z_i \chi(i) \leq \sqrt{ n \sum_{i=1}^n Z_i^2
(\chi(i))^2} \leq n^{1/2} B \sqrt{ \sum_{i=1}^n (\chi(i))^2} \leq
C n^{1/2},\eea for some $0<C<\infty$.
Hence we obtain (\ref{0922r}) in this case also. This proves part
(a).

Now we prove part (b). Suppose that for some $k \in \N$ $\chi$ is
$k$-supercritical. Again, we use the notation of (\ref{1022a}).
 By
(\ref{super}), we have that $\Var[S_n] \to \infty$ as $n \to
\infty$. Then, by Lemma \ref{itlog}, taking
$X_i=-(Z_i-\Exp[Z_i])$, we have that, with probability one,
\[ S_n < -( \Var[S_n])^{1/2} (3 \log \log (\Var[S_n]) )^{1/2} ,\]
for only finitely many $n$.
But $\chi(n) \to 0$ as $n \to \infty$, so with probability one there
exists a sequence $c_1,c_2,\ldots$
such that $c_n \to \infty$ as $n \to \infty$
and $\Var[S_n] < n/c_n$ for all $n$. Thus, with probability one,
\bea
\label{1811a}
S_n \geq -n^{1/2} c_n^{-1/2} (3 \log \log n)^{1/2},
\eea
 for all but finitely many $n$.
So, with probability one, for all but finitely many $n$, \bea
\label{0922ab} \sum_{i=1}^n Z_i \chi(i) \geq \Exp[Z_1]
\sum_{i=1}^n \chi(i) - n^{1/2} c_n^{-1/2} (3 \log \log n)^{1/2}.
\eea By (\ref{super}), we have that there exist $c \in (0,\infty)$
and $k \in \N$ such that for $n$ sufficiently large \bea
\label{0922abb} \sum_{i=1}^n \chi(i) \geq \frac{\sigma}{2
|\lambda|} \sum_{i=1}^n i^{-1/2} \varphi_k (i; c/2). \eea Then
from (\ref{0922abb}) with (\ref{0604b}) we obtain, for all $n$
sufficiently large \bea \label{0922aaby} \sum_{i=1}^n \chi(i)
 \geq \frac{\sigma}{|\lambda|} n^{1/2} \varphi_k(n;c/2) - \frac{3 \sigma }{|\lambda|} n^{1/2}.
\eea
Hence, with probability one, from (\ref{0922ab}) and (\ref{0922aaby}) we have that, for all
but finitely many $n$
\[
\sum_{i=1}^n Z_i \chi(i) \geq \frac{ \sigma \Exp[Z_1]}{ |\lambda|}  n^{1/2} \varphi_k (n; c/2)
- \frac{3 \sigma \Exp[Z_1]}{|\lambda|} n^{1/2} - n^{1/2} c_n^{-1/2} (3 \log \log n)^{1/2},\]
 which yields (\ref{0922t}). Thus we have proved part (b).

Finally, we prove part (c). Suppose now that $\Exp[Z_1]=0$. Again
use the notation of (\ref{1022a}). First, suppose that $\Var[S_n]
\leq C$ for all $n$, for some $0<C<\infty$. Then, we have that
(\ref{987a}) holds. On the other hand, suppose that $\Var[S_n] \to
\infty$ as $n \to \infty$. But, since $\chi(n) \to 0$ as $n \to
\infty$, we have that $\Var[S_n] = o(n)$. Applying Lemma
\ref{itlog} with $X_i = Z_i \chi(i)$ then yields (\ref{1007a}).
Thus the proof
of the lemma is complete.
$\square$ \\

\noindent {\bf Proof of Theorem \ref{cor3}.} First we examine the
recurrence and transience criteria for $\eta_t(\omega)$. For the
recurrent cases, we proceed in the second part of the proof to
analyse the stationary measure given in Lemma \ref{stat}, in order
to distinguish between null-recurrence and ergodicity (positive
recurrence). We work for a fixed environment $\omega$, that is, a
given realization of $p_i$ and $q_i$ for $i=1,2,\ldots$, as given
by (\ref{1006b}).

We aim to apply Lemmas \ref{crit1} and \ref{crit2}, and so we
construct a Lyapunov function $f$, that is, a function $f: \Z^+
\to \R^+$ such that $f(\eta_t(\omega))$ behaves as a martingale
(with respect to the natural filtration) for $\eta_t(\omega) \neq
0$. To do this, we proceed as follows.

For a given environment $\omega$, set $\Delta_1:=1$ and for
$i=2,3,\ldots$ let \bea \label{0818a}
 \Delta_i := \prod_{j=1}^{i-1} (p_j/q_j) = \exp
\sum_{j=1}^{i-1} \log{ (p_j/q_j)} , \eea and then set $f(0):=0$
and for $n=1,2,3,\ldots$ let \bea \label{0818b}
 f(n) := \sum_{i=1}^n \Delta_i .
\eea Note that $f(n) \geq 0$. Then, for fixed $\omega$, for $t \in
\Z^+$ and $n=1,2,\ldots$,
 \bean \Exp [ f(
\eta_{t+1}(\omega)) - f(\eta_t(\omega)) |
\eta_t(\omega) = n ] & = & p_n f(n-1) + q_n f(n+1) - f(n) \\
& = & q_n \Delta_{n+1} - p_n \Delta_n = 0 , \eean
i.e.~$f(\eta_t(\omega))$ is a martingale over $1,2,3,\ldots$.

We need to examine the behaviour of $f(n)$ as $n \to \infty$,
in order to apply Lemmas \ref{crit1} and \ref{crit2}. Recall from
(\ref{1006c}) that there exists
$n_0 \in \N$ such that for any $j>n_0$ and almost every realization of the
random environment $\omega$,
$p_j = \xi_j + Y_j \chi(j)$
and $q_j = 1- \xi_j -Y_j \chi(j)$.
Then, for $j$ sufficiently large, and a.e.~$\omega$
 \[\log p_j  =  \log ( \xi_j + Y_j \chi(j))  =  \log ( \xi_j ) + \xi_j^{-1} Y_j \chi(j)
 + O \left( (\chi(j))^2 \right) , \]
 and
\[
 \log q_j = \log ( 1-\xi_j - Y_j \chi(j))  =
\log ( 1-\xi_j ) -  (1-\xi_j)^{-1} Y_j \chi(j)
 + O \left( ( \chi(j))^2 \right) , \] so that for $j$ sufficiently large and a.e.~$\omega$
\bea
\label{0531b}
 \log (p_j /q_j) = \log \left( \frac{\xi_j}{1-\xi_j} \right) +
 \frac{Y_j}{\xi_j(1-\xi_j)} \chi(j) +O \left( (\chi(j))^2 \right)  .
\eea
Note that $\bbE [ \log (p_n/q_n) ] = O( \chi(n)) \to 0$ as $n \to \infty$, so that in this
sense
we asymptotically approach Sinai's regime.

Recall from (\ref{0520b}) that for $i=1,2,\ldots$,
$\zeta_i=\log(\xi_i/(1-\xi_i))$. From (\ref{0818b}), (\ref{0818a})
and (\ref{0531b}) we have, for $n$ sufficiently large, for
a.e.~$\omega$ \bea \label{0902a} f(n) = \sum_{i=1}^n \exp
\sum_{j=1}^{i-1} \left[
 \zeta_j  +
 \frac{Y_j}{\xi_j(1-\xi_j)} \chi(j) +O( (\chi(j))^2 )
\right]. \eea Note that for what follows the $O((\chi(j))^2)$
terms in (\ref{0902a}) can be ignored, since, when $\lambda \neq
0$ (where $\lambda$ is given by (\ref{0530a})), the other two
terms are dominant. Thus we need to examine the behaviour of the
two terms $\sum_{i=1}^n \zeta_i$ and $\sum_{i=1}^n
\frac{Y_i}{\xi_i(1-\xi_i)} \chi(i)$. This behaviour depends upon
the sign of $\lambda$, and the magnitude of the perturbation
$\chi$.

First suppose that for some $k \in \N$ $\chi$ is $k$-subcritical
(\ref{sub}). In this case, we show that in (\ref{0902a}) the term
involving the $\zeta_j$ is essentially dominant. We can apply
Lemma \ref{0530b} with $Z_i = Y_i \xi_i^{-1} (1-\xi_i)^{-1}$ (if
$\lambda>0$) or $Z_i = -Y_i \xi_i^{-1}(1-\xi_i)^{-1}$ (if
$\lambda<0$), and the boundedness property (\ref{bndd}), so that
(\ref{0922r}) implies that, for any $\eps>0$, for all but finitely
many $n$, for a.e.~$\omega$ \bea \label{0922xb} -n^{\eps} \leq
\sign( \lambda) \sum_{i=1}^n \frac{Y_i}{\xi_i (1-\xi_i)} \chi(i)
\leq \sigma n^{1/2} \varphi_k (n;-c/3), \eea with $c \in
(0,\infty)$ as given in (\ref{sub}). Also, from the Law of the
Iterated Logarithm (Lemma \ref{itlog}), we have that, for
a.e.~$\omega$, there are infinitely many values of $n$ for which
\bea \label{0922xa}
 \sum_{i=1}^n \zeta_i \geq \sigma n^{1/2} \varphi_k(n;-c/4).
\eea So from (\ref{0922xb}) and (\ref{0922xa}), we have that, for
a.e.~$\omega$, there are infinitely many values of $n$ such that,
if $\lambda>0$, \[ \sum_{i=1}^n \zeta_i + \sum_{i=1}^n
\frac{Y_i}{\xi_i (1-\xi_i)} \chi(i) \geq \sigma n^{1/2} \varphi_k
(n;-c/4) -n^\eps,\] and if $\lambda<0$,
\[ \sum_{i=1}^n \zeta_i + \sum_{i=1}^n
\frac{Y_i}{\xi_i (1-\xi_i)} \chi(i) \geq \sigma n^{1/2} (
\varphi_k (n;-c/4) -\varphi_k(n; -c/3)).\] Thus, by choosing
$\eps$ to be small, we have that for a.e.~$\omega$, there are
infinitely many values of $n$ such that
 \bea \label{0921e}
 \sum_{i=1}^n \zeta_i + \sum_{i=1}^n \frac{Y_i}{\xi_i (1-\xi_i)} \chi(i) \geq C n^{1/2}, \eea
for some $C$ with $0<C<\infty$. Thus from (\ref{0921e}),
(\ref{0818a}), and (\ref{0531b}), there are, for a.e.~$\omega$,
infinitely many values of $n$ for which $\Delta_n>1$, and hence as
$n \to \infty$ $f(n) \to +\infty$ for a.e.~$\omega$. Thus, by
Lemma \ref{crit2}, $\eta_t(\omega)$ is recurrent for
a.e.~$\omega$.

Now suppose that for some $k \in \N$ $\chi$ is $k$-supercritical
(\ref{super}). In this case, we show that the term in
(\ref{0902a}) involving $Y_j \xi_j^{-1} (1-\xi_j)^{-1}$ is
essentially dominant, and thus the sign of $\lambda$ determines
the behaviour. This time, from Lemma \ref{0530b} with $Z_i = Y_i
\xi_i^{-1} (1-\xi_i)^{-1}$ (if $\lambda>0$) or $Z_i = -Y_i
\xi_i^{-1} (1-\xi_i)^{-1}$ (if $\lambda<0$), and the boundedness
property (\ref{bndd}), we have that (\ref{0922t}) implies that,
for a.e.~$\omega$, for all but finitely many $n$, \bea
\label{0922xc} \sign(\lambda) \sum_{i=1}^n \frac{Y_i}{\xi_i
(1-\xi_i)} \chi(i) \geq \sigma n^{1/2} \varphi_k (n;c/3). \eea
Also, from the Law of
the Iterated Logarithm (Lemma \ref{itlog}), we have that, for
a.e.~$\omega$, there are only finitely many $n$ such that \bea
\label{0922xd}
 \sum_{i=1}^n \zeta_i \geq \sigma n^{1/2} \varphi_k(n;c/4).
\eea If $\lambda<0$, from (\ref{0922xc}) and (\ref{0922xd}), we
have that, for a.e.~$\omega$, there are only finitely many $n$
such that \bea \label{0921ee}
 \sum_{i=1}^n \zeta_i
+ \sum_{i=1}^n \frac{Y_i}{\xi_i (1-\xi_i)} \chi(i)  \geq \sigma
n^{1/2} \left( \varphi_k(n;c/4) -\varphi_k(n;c/3) \right) . \eea
So if $\lambda<0$, from (\ref{0921ee}), (\ref{0818a}), and
(\ref{0531b}), we have that for a.e.~$\omega$ there are only
finitely many values of $n$ for which
\[ \Delta_n \geq \exp \left(  - C_1 n^{1/2}  \right),\]
for some $C_1$, not depending on $\omega$, with $0<C_1<\infty$.
Thus for a.e.~$\omega$ there exists a constant $C_2$ (depending on
$\omega$) with $0<C_2<\infty$ such that
\[ f(n) \leq C_2 + \sum_{i=1}^\infty \exp \left(  - C_1 i^{1/2} \right),\]
which is bounded. So in this case, by Lemma \ref{crit1}, we have
that, for a.e.~$\omega$, $\eta_t(\omega)$ is transient.

On the other hand, if $\lambda>0$
then Lemma \ref{itlog} with (\ref{0922xc})
implies that for a.e.~$\omega$ there are infinitely many values of $n$ for which
\bea
\label{0921ek}
 \sum_{i=1}^n \zeta_i
+ \sum_{i=1}^n \frac{Y_i}{\xi_i (1-\xi_i)} \chi(i) \geq \sigma
n^{1/2} \left( \varphi_k(n;c/3) - \varphi_k(n;c/4)  \right) \geq
C_1 n^{1/2} , \eea for some $C_1$, not depending on $\omega$, with
$0<C_1<\infty$. So if $\lambda>0$, from (\ref{0921ek}),
(\ref{0818a}), and (\ref{0531b}) for a.e.~$\omega$ there are
infinitely many values of $n$ for which
\[ \Delta_n \geq \exp \left(  C_1 n^{1/2}  \right).\]
Thus $f(n) \to +\infty$ $\bbP$-a.s., and
 in this case we have that, for a.e.~$\omega$, $\eta_t(\omega)$ is recurrent, by
Lemma \ref{crit2}.

We now classify the recurrent cases further into ergodic (positive
recurrent) and null-recurrent. To determine ergodicity, we apply
Lemma \ref{stat}. Given $\omega$, and with $D(\omega)$ as defined
at (\ref{0520a}), we have
\[ D(\omega) = \sum_{i=1}^\infty \frac{1}{q_i} \exp
\left( -\sum_{j=1}^i \log{ (p_i/q_i) } \right) = \sum_{i=1}^\infty
\frac{1}{\Delta_{i+1} q_i} ,
\]
where $\Delta_i$ is as defined at (\ref{0818a}). By a similar
argument to (\ref{0531b}), we have that for $n$ sufficiently
large, for a.e.~$\omega$ \bean \frac{1}{\Delta_n} = \exp \left( -
\sum_{i=1}^{n-1} \zeta_i - \sum_{i=1}^{n-1} \frac{  Y_i}{\xi_i
(1-\xi_i)} \chi(i) + O \left(\sum_{i=1}^{n-1} (\chi(i))^2 \right)
\right). \eean We use similar arguments as in the proof of
recurrence and transience to analyse $D(\omega)$. First suppose
that for some $k \in \N$ $\chi$ is $k$-subcritical. Then, by a
similar argument to (\ref{0921e}), we have that for a.e.~$\omega$
there are infinitely many values of $i$ for which
\[ - \sum_{i=1}^n \zeta_i - \sum_{i=1}^n \frac{Y_i}{\xi_i (1-\xi_i)} \chi(i) \geq Cn^{1/2},\]
for $0<C<\infty$. Thus for a.e.~$\omega$ there are infinitely many
values of $n$ for which $(1/\Delta_{n+1})
>1$ and $(1/(\Delta_{n+1}q_n))>1$. Hence $D(\omega) = +\infty$ for
a.e.~$\omega$. So, for a.e.~$\omega$, by Lemma \ref{stat},
$\eta_t(\omega)$ is not ergodic.

Now suppose that for some $k \in \N$ $\chi$ is $k$-supercritical.
If $\lambda >0$, using similar arguments to before,
 we have that for a.e.~$\omega$ there are only finitely many $n$ for which
\bean - \sum_{i=1}^n \zeta_i -\sum_{i=1}^n \frac{Y_i}{\xi_i
(1-\xi_i)} \chi(i)  \geq \sigma n^{1/2} \left( \varphi_k(n;c/4)
-\varphi_k(n;c/3) \right) . \eean  So for a.e.~$\omega$ there are
only finitely many values of $n$ for which
\[ (1/\Delta_n) \geq \exp \left(  - C_1 n^{1/2}  \right),\]
for some $0<C_1<\infty$. Thus for a.e.~$\omega$ there exists a
constant $C_2$ (depending on $\omega$) with $0<C_2<\infty$ such
that
\[ D(\omega) \leq C_2 + \sum_{i=1}^\infty \exp \left(  - C_1 i^{1/2}  \right),\]
which is bounded. So in this case, for a.e.~$\omega$, by Lemma
\ref{stat}, $\eta_t(\omega)$ is ergodic.
This completes the proof of Theorem \ref{cor3}. $\square$ \\

\noindent {\bf Proof of Theorem \ref{cor2}.}  First we prove parts
(i) and (iii). Suppose that, for all $n$ sufficiently large,
$\lambda \chi(n) \geq c n^{-1/2}(\log \log n)^{1/2}$, for some
$c>c_{\rm crit}$ where $c_{\rm crit} = \sigma 2^{-1/2}$. Then we
see that $\chi$ is $k$-supercritical (\ref{super}) for
$k=2,3,\ldots$, since, for example \bean
 \frac{c}{|\lambda|} n^{-1/2}(\log \log n)^{1/2} & = &
\frac{c}{c_{\rm crit}} \frac{\sigma}{2 |\lambda|} n^{-1/2} ( 2 \log \log {n} )^{1/2} \\
& \geq & \frac{\sigma}{2 |\lambda|} n^{-1/2} (2 \log \log {n} + 4 \log \log \log n)^{1/2} ,\eean
for $n$ sufficiently large and $c > c_{\rm crit}$. Hence (i) follows from part (i) of Theorem \ref{cor3}.
Similarly, (iii) follows from part (iii) of Theorem \ref{cor3}.

For part (ii), suppose that $|\lambda| \chi(n) \leq c n^{-1/2}(\log \log n)^{1/2}$
for all $n$ sufficiently large, $c \leq c_{\rm crit}$.
Then we see that $\chi$
is $k$-subcritical (\ref{sub}) for $k=2,3,\ldots$, since, for example
\bean
 \frac{c}{|\lambda|} n^{-1/2}(\log \log n)^{1/2} & \leq &
\frac{\sigma}{2 |\lambda|} n^{-1/2} ( 2 \log \log {n} )^{1/2} \\
& \leq & \frac{\sigma}{2 |\lambda|} n^{-1/2} (2 \log \log {n} + 2 \log \log \log n)^{1/2} ,\eean
for $n$ sufficiently large. Then part (ii) of Theorem \ref{cor3} gives part (ii) of Theorem \ref{cor2},
and the proof
of Theorem \ref{cor2} is complete. $\square$ \\

\noindent {\bf Proof of Theorem \ref{cor0}.}  From Lemma
\ref{itlog}, we have that for a.e.~$\omega$ there are infinitely
many values of $n$ for which \bea \label{664} \sum_{i=1}^n \zeta_i
\geq \sigma n^{1/2} (\log \log n)^{1/2}.\eea By a similar
argument to (\ref{0531b}), but keeping track of higher order terms
in the Taylor series, we have that now \bea \label{1102h}
\log(p_i/q_i) = \zeta_i + \sum_{r=1}^\infty \frac{1}{r} Y_i^{r}
\left( \frac{1}{(1-\xi_i)^r}+\frac{(-1)^{r+1}}{\xi_i^r} \right)
(\chi(i))^{r} . \eea By the condition $Y_1/\xi_1 \eqd
-Y_1/(1-\xi_1)$, we have that the expectation of the sum on the
right of (\ref{1102h}) is zero. Hence we can apply part (c) of
Lemma \ref{0530b} with \bea \label{0623a}
 Z_i = \sum_{r=1}^\infty \frac{1}{r} Y_i^{r} \left( \frac{1}{(1-\xi_i)^r}+\frac{(-1)^{r+1}}{\xi_i^r}
\right) (\chi(i))^{r-1} \eea to obtain that for all but finitely
many $n$, for a.e.~$\omega$ \bea \label{665} \sum_{i=1}^n
\sum_{r=1}^\infty \frac{1}{r} Y_i^{r} \left(
\frac{1}{(1-\xi_i)^r}+\frac{(-1)^{r+1}}{\xi_i^r} \right)
(\chi(i))^{r} \geq -\eps n^{1/2} ( \log \log n)^{1/2},\eea and
by choosing $\eps$ sufficiently small we have from (\ref{1102h}),
(\ref{664}) and (\ref{665}) that, for a.e.~$\omega$ there are
infinitely many values of $n$ for which
\[ \sum_{i=1}^n \log (p_i/q_i) \geq C n^{1/2} (\log \log
n)^{1/2},\] for $0<C<\infty$.
Thus with $\Delta_n$ defined at
(\ref{0818a}), we have that for a.e.~$\omega$ there are infinitely
many values of $n$ for which
\[ \Delta_n \geq \exp \left( C n^{1/2} (\log \log n)^{1/2} \right) ,\]
and so $f(n) \to +\infty$ $\bbP$-a.s., and so, by Lemma
\ref{crit2}, $\eta_t(\omega)$ is recurrent for a.e.~$\omega$.

To prove null-recurrence, it remains to show that the Markov chain
is not ergodic. Consider $D(\omega)$ as defined at (\ref{0520a})
again. From Lemma \ref{itlog}, we have that for a.e.~$\omega$
there are infinitely many values of $n$ for which
\[ -\sum_{i=1}^n \zeta_i \geq \sigma n^{1/2} (\log \log n)^{1/2}.\]
From part (c) of Lemma \ref{0530b}
with $Z_i$ as at (\ref{0623a})
we have that
for all but finitely many $n$, for a.e.~$\omega$
\[ - \sum_{i=1}^n Z_i \chi(i) \geq -\eps n^{1/2} ( \log \log n)^{1/2},\]
and by choosing $\eps$ sufficiently small we have that for a.e.~$\omega$ there
are infinitely many values of $n$ for which
\[ (1/\Delta_n) \geq \exp \left( C n^{1/2} (\log \log n)^{1/2} \right) ,\]
for some $0<C<\infty$, and so $D(\omega)= +\infty$ $\bbP$-a.s.
Thus, by Lemma \ref{stat}, the Markov chain is $\bbP$-a.s.~not
ergodic. Thus, for a.e.~$\omega$, $\eta_t(\omega)$ is
null-recurrent. $\square$ \\

\noindent {\bf Proof of Theorem \ref{thm1022}.}  Parts (i) and
(ii) follow easily with the methods used in the proof of Theorem
\ref{cor3}. We prove part (iii). By a similar argument to
(\ref{0531b}), we have that now \bea \label{1102a} \log(p_i/q_i) =
\sum_{r=1}^\infty \frac{4^r}{2r-1} Y_i^{2r-1} (\chi(i))^{2r-1} = 4
Y_i \chi(i) +O((\chi(i))^3) . \eea Since $Y_1 \eqd -Y_1$, we have
that all odd powers of $Y_1$ have zero expectation, so that the
expectation of the right hand side of (\ref{1102a}) is zero. Thus
it is clear that for a.e.~$\omega$ there are infinitely many
values of $n$ for which $\sum_{i=1}^n \log(p_i/q_i) \geq 0$, and
hence $\Delta_n \geq 1$, and so $f(n) \to +\infty$ for
a.e.~$\omega$, and we have $\bbP$-a.s.~recurrence, by Lemma
\ref{crit2}.

To prove null recurrence, it remains to show that the Markov chain
is not ergodic. Once more, consider $D(\omega)$ as defined at
(\ref{0520a}). By a similar argument to above, for a.e.~$\omega$
there are infinitely many values of $n$ for which $\sum_{i=1}^n
\log(p_i/q_i) \leq 0$ and hence $(1/\Delta_n) \geq 1$, and so
$D(\omega)= +\infty$ for a.e.~$\omega$. Thus, by Lemma \ref{stat},
the Markov chain is $\bbP$-a.s.~not ergodic. This completes the
proof of part (iii).

We now prove part (iv). Once again we analyse the properties of the
expression (\ref{1102a}). 
Suppose that $\chi(n)= a n^{-\beta}$ for $a>0$, $\beta>0$. Now
suppose that $0<\beta<1$ and that $\bbE[Y_1]<0$. Then from
(\ref{1102a}), we have that there exist  $0<C_1<\infty$,
$0<C_2<\infty$ such that
\[ -C_1n^{1-\beta} \leq \bbE \sum_{i=1}^n \log(p_i/q_i) \leq -C_2n^{1-\beta}.\]
If $\beta > 1/2$, then, by the boundedness of $Y_1$, we have
\[ \sup_n \bbE \left| \sum_{i=1}^n \log(p_i/q_i) - \bbE  \sum_{i=1}^n \log(p_i/q_i)
\right|^k < \infty ,\]
for all $k \in \N$, so that $\bbP$-a.s.,
\[  \left| \sum_{i=1}^n \log(p_i/q_i) - \bbE  \sum_{i=1}^n \log(p_i/q_i)
\right| < n^\eps ,\]
for all but finitely many $n$, and any $\eps>0$. So, for all but finitely many $n$, for a.e.~$\omega$
\[ \Delta_n \leq \exp \left( -Cn^{1-\beta}+n^{\eps} \right),\]
for some $C$ with $0<C<\infty$,
so that, for $\eps$ small enough, $f(n) = \sum_{i=1}^n \Delta_i$
is bounded for a.e.~$\omega$, which implies that $\eta_t(\omega)$
is $\bbP$-a.s.~transient, by Lemma \ref{crit1}. 

Also, if $\beta=1/2$ we
have from (\ref{1102a}) that there exist
$0<C_1<\infty$, $0<C_2<\infty$ such that
\[ C_1 \log n \geq \Var \sum_{i=1}^n \log(p_i/q_i) \geq C_2 \log n \to \infty,\]
as $n \to \infty$, and then we can apply Lemma \ref{itlog} to obtain, for a.e.~$\omega$
\[ \sum_{i=1}^n \log(p_i/q_i) \leq -C_1 n^{1/2} + C_2 (\log n)^{1/2} (\log \log {n})^{1/2} ,\]
for constants $0<C_1<\infty$, $0<C_2<\infty$ (depending on $\omega$)
and all but finitely many $n$. Thus $f(n)$ is $\bbP$-a.s.~bounded, and so we have $\bbP$-a.s.~transience by Lemma \ref{crit1}. 

Finally, if 
 $0< \beta
< 1/2$, from (\ref{1102a}), we have that there exist
$0<C_1<\infty$, $0<C_2<\infty$ such that
\[ C_1 n^{1-2\beta} \geq \Var \sum_{i=1}^n \log(p_i/q_i) \geq C_2 n^{1-2\beta} \to \infty,\]
as $n \to \infty$, and then by Lemma \ref{itlog} we obtain, for a.e.~$\omega$
\[ \sum_{i=1}^n \log(p_i/q_i) \leq -C_1 n^{1-\beta} + C_2 n^{(1/2)-\beta} (\log \log {n})^{1/2} ,\]
for constants $0<C_1<\infty$, $0<C_2<\infty$ (depending on $\omega$)
and all but finitely many $n$. So once again we
have $f(n)$ is $\bbP$-a.s.~bounded, and so we have $\bbP$-a.s.~transience by Lemma \ref{crit1}. This proves part (c).

To prove part (a), we apply Lemma \ref{stat}.
Suppose that $\bbE[Y_1]>0$. By similar arguments to above, this time we have that for a.e.~$\omega$
\[ \sum_{i=1}^n \log(p_i/q_i) \geq Cn^{1-\beta},\]
for some $0<C<\infty$ and all but finitely many $n$. Thus, for a.e.~$\omega$, for all but finitely many $n$,
\[ \frac{1}{\Delta_n} = \exp \left( -\sum_{i=1}^{n-1} \log(p_i/q_i) \right)  \leq \exp \left( - Cn^{1-\beta} \right),\]
and so, for $D(\omega)$ as defined at (\ref{0520a}),
 $D(\omega) < \infty$ $\bbP$-a.s., and so, by
Lemma \ref{stat}, the Markov chain is $\bbP$-a.s.~ergodic,
proving part (a).

Finally, we prove part (b). Suppose that $\beta>1$. Now, since $-1
\leq Y_i \leq 1$ and $\chi(n)=an^{-\beta}$, we have from
(\ref{1102a}) that there exists a constant $C_1$ (not depending on
$\omega$) with $0<C_1<\infty$ such that, for a.e.~$\omega$,
\[ \left| \sum_{i=1}^n  \log(p_i/q_i) \right| \leq C_1 \sum_{i=1}^n i^{-\beta} \leq C_2,\]
for finite positive $C_2$, not depending on $\omega$ or $n$, this
last inequality following since $\beta>1$. Thus for a.e.~$\omega$,
for each $n$,
\[ 0<\exp(-C_2) \leq \exp \left( \sum_{i=1}^n \log(p_i/q_i) \right) \leq \exp (C_2)<\infty ,\]
so that for each $n$, $\Delta_n$ and $1/\Delta_n$ are each bounded
strictly away from $0$ and from $\infty$, so that
$\bbP$-a.s.~$f(n) \to + \infty$ as $n \to \infty$, and
$D(\omega)=+\infty$ $\bbP$-a.s. Thus by Lemma \ref{crit1} the
Markov chain is $\bbP$-a.s.~recurrent, and by Lemma \ref{stat}
$\bbP$-a.s.~not ergodic. Thus, for a.e.~$\omega$, $\eta_t(\omega)$
is null-recurrent. This completes the proof of Theorem
\ref{thm1022}. $\square$

\vskip 2mm
\begin{center}
{\bf Acknowledgements}
\end{center}
\vskip 1mm
 AW began this work while at the University of Durham, supported by an EPSRC
doctoral training account, and it was completed when he was at the University of Bath.

\end{document}